# Conway's Circle Theorem: A Short Proof, Enabling Generalization to Polygons

Eric Braude
*MET Computer Science Department*
*Boston University*
Boston, USA
orcid.org/0000-0002-1630-5509

*Abstract*— John Conway's Circle Theorem is a gem of plane geometry: the six points formed by continuing the sides of a triangle beyond every vertex by the length of its opposite side, are concyclic. The theorem has attracted several proofs, even adorned Mathcamp T-shirts. We present a short proof that views the extended sides as equal tangents of the incircle, a perspective that enables generalization to polygons.

## I. INTRODUCTION

Several proofs of Conway's Circle Theorem exist. A recent proof by distinguished geometer Doris Schattschneider [7] presents what she calls "a nicer, more convincing" proof than her first—which had already been praised for its "aesthetic appeal, … perhaps the one Euclid would have worked out had he noticed Conway's theorem himself" [7]. Schattschneider's second proof is based on elegant constructions joining points of the hexagon. In May 2020, Aperiodical published a proof without words by Colin Beveridge [1]. In [7], Baker calls this "one of the most beautiful Proofs Without Words I've ever seen". It constructs line segments joining the six points to the incentre. A video with words appeared in June 2020 [2]. Alex Ryba [8] produced a short proof appealing to simple constructions and overlapping isosceles triangles.

This paper is based on the following proof of Conway's Circle Theorem for triangle ABC with respective opposite side lengths a, b, and c, as shown in Fig. 1.

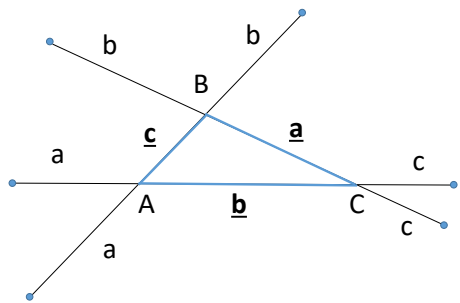

Fig. 1. Conway's Circle Theorem.

## II. PROOF OF CONWAY'S CIRCLE THEOREM

As shown in Fig. 2, let $I$ be the incenter of *ABC*, $P_{CA}$ the point on the extension of *CA* at a distance $a$ from *A*, and $\Omega$ the circle with center $I$ containing $P_{CA}$. Let $P_{CB}$, $P_{AB}$, $P_{AC}$, $P_{BC}$, and $P_{BA}$ be the intersections of $\Omega$ with the extensions of sides *CB*, *AB*, *AC*, *BC*, and *BA* respectively. Chords $P_{CA}P_{AC}$, $P_{AB}P_{BA}$, and $P_{BC}P_{CB}$ are rotations of each other around *I*, and are thus equal. Since each chord pair is symmetric about the diameter through its intersection, $AP_{BA} = a$, $CP_{BC} = c'$, and $BP_{AB} = b'$. The chord equality thus yields $b' + a + c' = b' + c + a = c' + b + a$. Thus, $c' = c$, $b' = b$, and Conway's theorem follows.

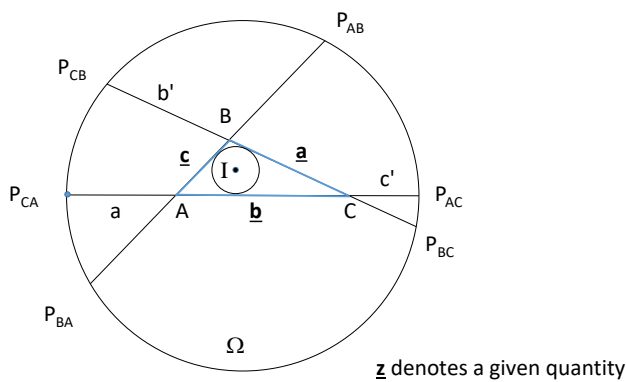

Fig. 2. Proof of Conway's Circle Theorem.

Generalizations of Conway's Circle Theorem have been investigated by several authors; for example, Capitan [3], who extended it to conics. The generalization in Theorem 1 below refers to the circles centered at the incenter and containing the incircle. We will call these Conway Circles. Theorem 1 reduces to Conway's Circle Theorem when $x_A = \lambda_A$, $x_B = \lambda_B$, and $x_C = \lambda_C$. In this case, the chord lengths equal the perimeter, an observation that will echo in subsequent theorems.

## III. THEOREM 1: CONWAY CIRCLES FOR TRIANGLES

Let *T* be a triangle with vertices *A*, *B*, and *C*, and respective opposite side lengths $\lambda_A$, $\lambda_B$, and $\lambda_C$. Suppose that the two sides passing through *A* (resp. *B*, *C*) are extended by a distance $x_A$ (resp. $x_B$, $x_C$), as illustrated in Fig. 3. The six endpoints of these extensions lie on a common circle centered at the incenter of *T* if and only if $x_U - x_V = \lambda_U - \lambda_V$ for all vertices *U* and *V* of *T*.

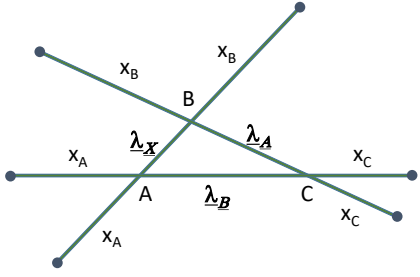

Fig. 3. Conway circles for triangles.

## IV. Proof of Theorem 1: Conway Circles for triangles

Let $I$ be $T$'s incenter, $A'$ the point on the extension of $CA$ defined by $x_A$, and $\Omega$ the circle centered at $I$ that contains $A'$. For each pair of vertices $U$ and $V$, let $y_{UV}$ be the segment length on the extension of $UV$ defined by $V$ and $\Omega$, as shown in Fig. 4.

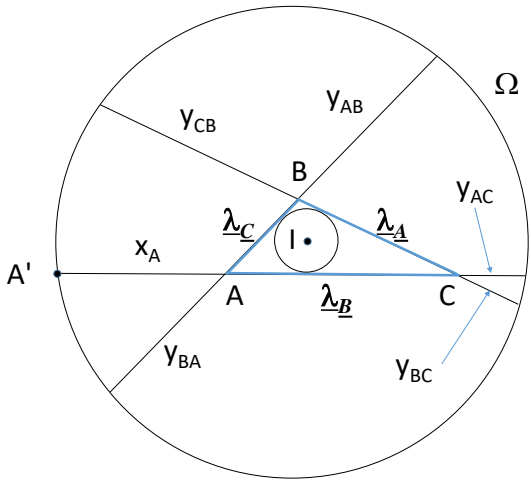

Fig. 4. Notation for proof of Theorem 1.

Since the three chords are tangential to the incircle, they are rotations of each other about $I$, and thus of equal length. Since each pair of these chords is symmetrical about the diameter through its intersection, $y_{BA} = x_A$, $y_{BC} = y_{AC}$, and $y_{AB} = y_{CB}$. We can thus label $y_{BC}$ and $y_{AC}$ as $y_C$; also, $y_{AB}$ and $y_{CB}$ as $y_B$. Thus,

$$y_B + \lambda_A + y_C = x_A + \lambda_B + y_C$$

and so

$$y_B = x_A + \lambda_B - \lambda_A$$

Similarly,

$$y_U = x_V + \lambda_U - \lambda_V \text{ for all vertices } U \text{ and } V \text{ of } T \quad (1)$$

To prove sufficiency: given $x_A$, $x_B$, and $x_C$ satisfying $x_V - x_W = \lambda_V - \lambda_W$ for all vertices $V$ and $\Omega$ of $T$, we have, using (1),

$$x_B = x_A + \lambda_B - \lambda_A = y_B$$

Similarly, $x_C = y_C$, so $\Omega$ is a circle centered at $I$, passing through the six points illustrated in Fig. 3.

Conversely, if there is a circle with incenter $I$ intersecting the six points specified by $x_A$, $x_B$, and $x_C$, then it must be $\Omega$ because the latter is the unique circle centered at the incenter and intersecting at $A'$, which is defined by $x_A$. Consequently, $x_Z = y_Z$ for $Z = A, B$, and C; and $x_U - x_V = y_U - y_V = \lambda_U - \lambda_V$ for all vertices $V$ and $\Omega$ of $T$ by (1).

The circle $\Omega$ in Theorem 1 reduces to the incircle when

$$(-x_A) + (-(\lambda_B - \lambda_A + x_A)) = \lambda_C, \text{ i.e., } x_A = \tfrac{1}{2}(\lambda_A - \lambda_B - \lambda_C).$$

For every $x_A < \tfrac{1}{2}(\lambda_A - \lambda_B - \lambda_C)$, the resulting circle coincides with one for which $x_A > \tfrac{1}{2}(\lambda_A - \lambda_B - \lambda_C)$ and so we can assume the latter. Since the results concern *all* circles at the incenter at least as large as the incircle, we will refer to them as "Conway circles."

Theorem 1 provides corollaries by selecting interesting values of $x_V$. For example, Corollary 1 follows by taking $x_A = 0$, and is illustrated in Fig. 5 for $\lambda_B > \lambda_A$ and $\lambda_C < \lambda_A$.

## V. Corollary 1: Conway Circles at Triangle Vertices

Given a triangle with vertices $A$, $B$, and $C$, and respective opposite side lengths $\lambda_A$, $\lambda_B$, and $\lambda_C$, if sides are extended by $\lambda_B - \lambda_A$ at $B$ and by $\lambda_C - \lambda_A$ at $C$, then the resulting five points, including $A$, are concyclic, with center at the incenter.

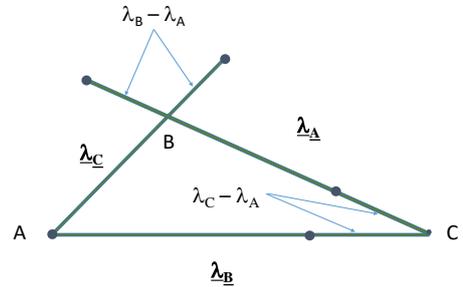

Fig. 5. Conway circles at a vertex.

A *tangential* polygon is one for which a circle exists whose sides are tangential to it. Every convex polygon $P$ corresponds to a set of tangential polygons with sides that are pairwise parallel to $P$. (On any circle $C$, construct successive tangents to $C$, each parallel to successive sides of $P$.) We will use the chord symmetry argument described above for tangential polygons with an odd number of sides ("odd polygons"). We then use this result to generalize for even tangential polygons. Theorem 2 generalizes Theorem 1. Its proof proceeds like that for Theorem 1, progressing around the polygon once for odd indices and then, because $n$ is odd, continuing in the same manner for even indices. This process uses a modified *mod* function $\mu()$.

## VI. THEOREM 2: CONWAY CIRCLES FOR ODD TANGENTIAL POLYGONS

Let $P$ be a tangential polygon with an odd number $n$ of vertices $V_1, V_2, \ldots, V_n$, define $\lambda_i$ as $V_iV_{i+1}$ for $i = 1, 2, \ldots, n-1$ and $\lambda_n$ as $V_nV_1$, illustrated in Fig. 6. Suppose that each pair of sides ending at $V_i$ are extended by length $x_i$. The $2n$ endpoints of these extensions lie on a common circle centered at the incenter of $P$ if and only if

$$x_i - x_{\mu(i+2,\, n)} = \lambda_{\mu(i+1,\, n)} - \lambda_i \text{ for } 0 < i \le n \quad (2)$$

—where $\mu(z, m)$ is defined as $z \bmod m$ for $z \ne m$, and $m$ otherwise.

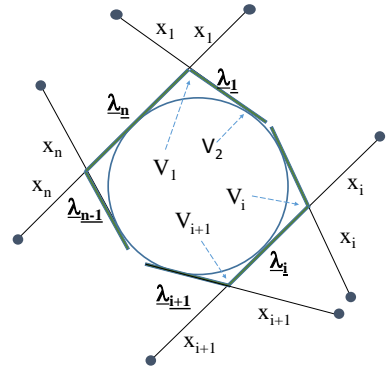

**z** denotes a given quantity

Fig. 6. Conway circles for odd tangential polygons.

## VII. PROOF OF THEOREM 2: CONWAY CIRCLES FOR ODD TANGENTIAL POLYGONS

As illustrated in Fig. 7, let $I$ be the incenter of the given polygon, $V_1'$ the point on the extension of $\lambda_1$ from $V_1$ defined by $x_1$, and $\Omega$ the circle, centered at $I$, containing $V_1'$. For $1 \le i \le n$, let $y_i$ and $y_i'$ be the segment lengths on the side extensions through $V_i$ defined by their intersections with $\Omega$.

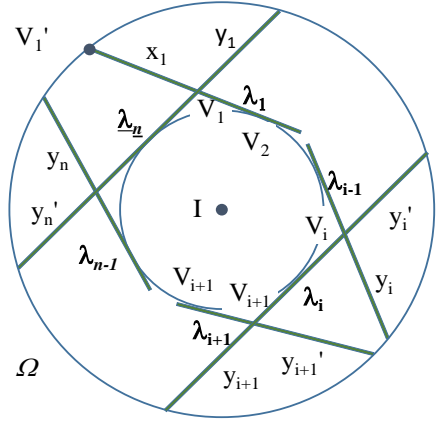

Fig. 7. For proof of Conway Circles for odd tangential polygons

As in the proof of Theorem 1, $y_1 = x_1$, $y_i = y_i'$ for all $2 \le i \le n$, and the $n$ chords are equal. Thus, for $1 \le i \le n-2$,

$$y_i + \lambda_i + y_{i+1} = y_{i+1} + \lambda_{i+1} + y_{i+2}$$

and so $y_i - y_{i+2} = \lambda_{i+1} - \lambda_i$

For $1 \le i < n-2$ and $i = n$, this establishes

$$y_i - y_{\mu(i+2,\, n)} = \lambda_{\mu(i+1,\, n)} - \lambda_i \quad (3)$$

Equation (3) is also satisfied for $i = n-1$ and $i = n-2$, since it reduces to, respectively,

$$y_{n-1} - y_{\mu(n+1,\, n)} = \lambda_{\mu(n,\, n)} - \lambda_{n-1}$$

and $\quad y_{n-2} - y_{\mu(n,\, n)} = \lambda_{\mu(n-1,\, n)} - \lambda_{n-2}$

i.e., $\quad y_{n-1} - y_1 = \lambda_n - \lambda_{n-1}$

and $\quad y_{n-2} - y_n = \lambda_{n-1} - \lambda_{n-2}$

But these follow from the chord length equalities

$y_{n-1} + \lambda_{n-1} + y_n = y_n + \lambda_n + y_1$ and $y_{n-2} + \lambda_{n-2} + y_{n-1} = y_{n-1} + \lambda_{n-1} + y_n$ resp.

To prove *sufficiency* in Theorem 2, assume that there is a circle $V$, centered at the incenter, that intersects the sides extended from $V_i$ at distances $x_i$ for $1 \le i \le n$. Thus, $V=W$ because $W$ is defined by the incenter and the distance $x_1$. The relationships (3) follow as above for Theorem 1, which concludes the sufficiency.

To prove the *necessity*, assume equation (2) for some $x_1, x_2, \ldots, x_n$. We thus have

$$x_1 - x_{\mu(3,\, n)} = \lambda_{\mu(2,\, n)} - \lambda_i$$

Constructing $W$ from $x_1$ as above, the consequence $x_1 = y_1$, and equation (3) together imply

$$\lambda_{\mu(2,\, n)} - \lambda_i = x_1 - y_{\mu(3,\, n)}$$

so $x_3 = y_3$. Similarly, $x_5 = y_5, x_7 = y_7, \ldots$, and $x_n = y_n$. For $i = n$, equation (2) becomes

$$x_n - x_{\mu(n+2,\, n)} = \lambda_{\mu(n+1,\, n)} - \lambda_n$$

i.e., since $x_n = y_n$, we have $y_n - x_2 = \lambda_1 - \lambda_n$.
But the equality of the chords in Fig. 7 implies

$$y_n + \lambda_n + y_1 = y_1 + \lambda_1 + y_2$$

so we have $y_n - x_2 = y_n - y_2$. Thus, $x_2 = y_2$, and the corresponding equalities continue for the even indices.

Corollary 2 below is easily recognized as a direct generalization of Conway's Circle Theorem for triangles, the chord length equaling the perimeter of the polygon. It follows from Theorem 2 by verifying equation (2).

## VIII. COROLLARY 2: CONWAY'S CIRCLE FOR ODD POLYGONS

Let $P$ be a tangential polygon with vertices $V_1, V_2, \ldots, V_n$, $n$ odd, and $\lambda_i$ denoting $V_iV_{i+1}$ for $i = 1, 2, 3, \ldots, n$, as shown in Fig. 8. Suppose that each pair of sides ending at $V_i$ are extended by length $\Sigma\{\lambda_i: (1 \leq i < k \wedge i \text{ odd}) \vee (k < i \leq n \wedge i \text{ even})\}$ for odd $k$, and $\Sigma\{\lambda_i: (1 \leq i < k \wedge i \text{ even}) \vee (k < i \leq n \wedge i \text{ odd})\}$ for even $k$. Then the $2n$ endpoints of these extensions lie on a common circle centered at the incenter of $P$.

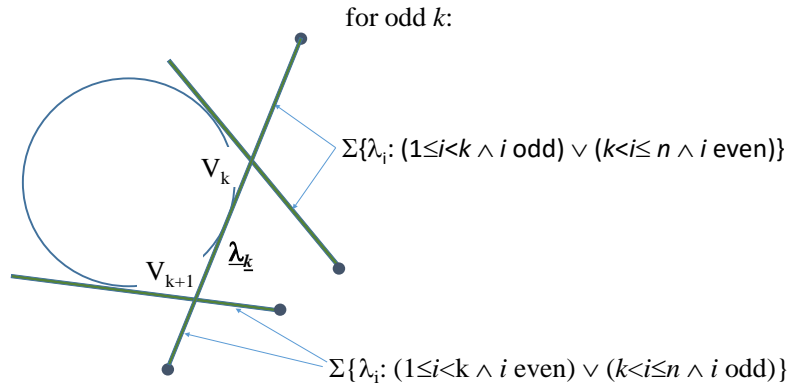

Fig. 8. Conway's circle for odd tangential polygons.

Corollary 3 is the application of Corollary 2 to pentagons, in which the generalization of Conway's Circle Theorem is quite graphic. It follows by taking $(\lambda_1, \lambda_2, \lambda_3, \lambda_4, \lambda_5) = (a, b, c, d, e)$.

## IX. COROLLARY 3: CONWAY'S CIRCLE FOR TANGENTIAL PENTAGONS

Let *ABCDE* be a tangential pentagon with opposite side lengths *a*, *b*, *c*, *d*, and *e* respectively. If the sides are extended at *A*, *B*, *C*, D, and *E* by $b + e$, $a + c$, $b + d$, $e + c$, and $a + d$ respectively, then the 10 resulting points are concyclic, with center at the pentagon's incenter. Fig. 9 illustrates this.

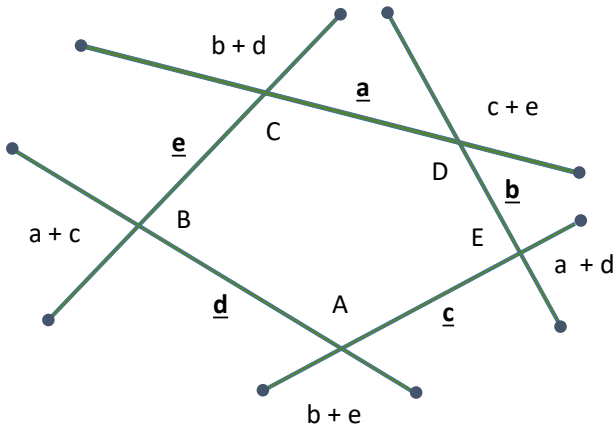

Fig. 9. Conway's Circles for tangential pentagons.

We will produce Conway circles next for *even* tangential polygons. Given an even sequence of lengths that form a tangential polygon, there are infinitely many such tangential polygons, as shown in [6]. The Conway circle formulae are slightly more complicated as a result. We will confine our result to necessary conditions, which continue to generalize the Conway Circle Theorem for triangles.

## X. THEOREM 3: CONWAY CIRCLES FOR EVEN TANGENTIAL POLYGONS

Given an even-sided tangential polygon $V_1, V_2, …, V_n$, with $h_1$ defined as $V_1V_2$, $h_2$ as $V_2V_3$, …, and $h_m$ as $V_mV_1$, in which the incircle is incident on $V_mV_1$ at a distance $h_0$ from $V_1$, extend each side $V_mV_{m+1}$, with odd $m$, by quantity (4) below at $V_m$ and by quantity (5) at $V_{m+1}$, then the $2n$ points so formed are concyclic, with center at the polygon's incenter. This is illustrated in Fig. 10.

$$\Sigma\{h_i: (1 \leq i < k \wedge i \text{ odd}) \vee (k < i \leq m \wedge i \text{ even})\} - h_0 \quad (4)$$

$$\Sigma\{h_i: (1 \leq i < k \wedge i \text{ even}) \vee (k < i \leq m \wedge i \text{ odd})\} + h_0 \quad (5)$$

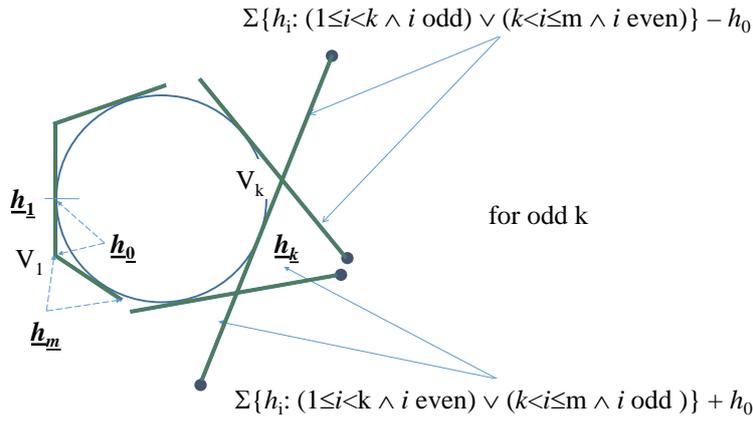

Fig. 10. Conway's circle for even tangential polygons.

## XI. PROOF OF THEOREM 3: CONWAY CIRCLES FOR EVEN TANGENTIAL POLYGONS

As shown in Fig. 11, we introduce a distance $h_0'$ from $V_1$ on $V_1V_2$ a little further than $h_0$ ($h_0'$ will converge to $h_0$), defining the point $V_2'$, then continuing the existing tangential polygon from $V_2'$ instead of $V_2$. This replaces the segments $V_1V_2$ and $V_2V_3$ with segments $V_1V_2'$, $V_2'V_2''$, and $V_2''V_3$. These have lengths $h_0'$, $h_1'$, and $h_2'$ respectively, say.

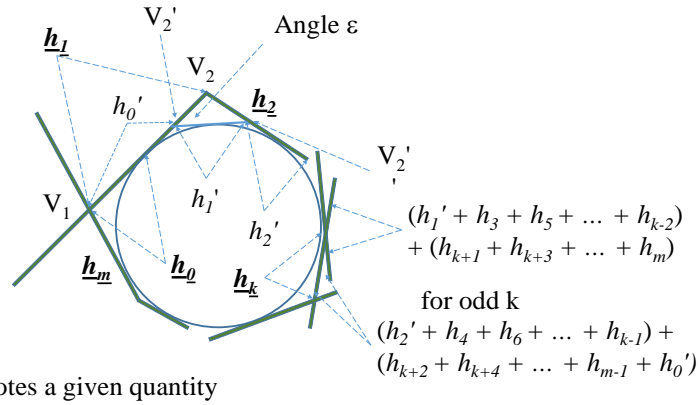

Fig. 11. Proving Conway's circle for even tangential polygons.

The resulting polygon has an odd number of sides, so we can apply Corollary 2 with $h_0'$ replacing $\lambda_n$ of Theorem 2, $h_1'$ replacing $\lambda_1$, $h_2'$ replacing $\lambda_2$, and $h_i$ replacing $\lambda_i$ for the remaining sides. The expressions shown in Figure 11 result. As angle $\varepsilon \to 0$, we have $h_0' \to h_0$, $h_1' \to h_1 - h_0$, $h_2' \to h_2$, and the following limits hold for the two expressions displayed in Fig. 11, proving Theorem 3:

$(h_1' + h_3 + h_5 + \ldots + h_{k-2}) + (h_{k+1} + h_{k+3} + \ldots + h_m)$
$\to \Sigma\{h_i: (1 \leq i < k \land i \text{ odd}) \lor (k < i \leq m \land i \text{ even})\} - h_0$

and

$(h_2' + h_4 + h_6 + \ldots + h_{k-1}) + (h_{k+2} + h_{k+4} + \ldots + h_{m-1} + h_0')$
$\to \Sigma\{h_i: (1 \leq i < k \land i \text{ even}) \lor (k < i \leq m \land i \text{ odd})\} + h_0$

Corollary 4 follows by taking $n = 4$, and $h_1$, $h_2$, $h_3$, and $h_4$ = $a$, $b$, $c$, and $d$ respectively in Theorem 3.

## XII. COROLLARY 4: CONWAY CIRCLE FOR TANGENTIAL QUADRILATERALS

For any tangential quadrilateral with side lengths $a$, $b$, and $c$, and the tangency on the remaining side—with length $d$—at a distance $d_0$ from a vertex, extend the sides between $a$ and $b$, $b$ and $c$, $c$ and $d$, and $d$ and $a$, by $c + d - d_0$, $a + d_0$, $b + d - d_0$, $b + d - d_0$, and $b + d_0$ respectively. Then the resulting eight points, as illustrated in Fig. 12, are concyclic.

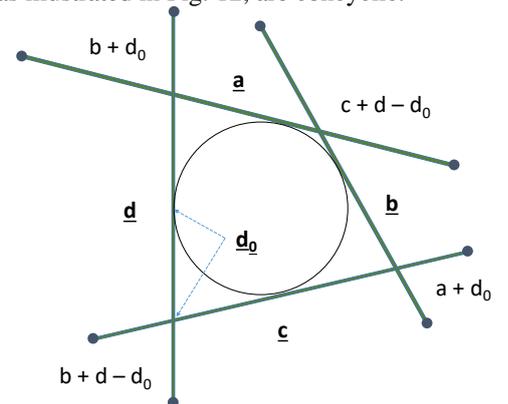

Fig. 12. Conway circles for a tangential quadrilateral.

In summary, when we view the side extensions in Conway's pretty circle theorem as equal chords tangential to the incircle, a clear perspective emerges, generalizable to tangential polygons

## Acknowledgment

The author is grateful to Matt Baker for his encouragement and for his elegant reformulation of Theorem 1.
## REFERENCES

[1] Beveridge, C. Conway's Circle, a proof without words, retrieved December 8, 2020 from https://aperiodical.com/2020/05/the-big-lock-down-math-off-match-14

[2] Beveridge, C. Conway's Circle Theorem: a proof, this time with words, retrieved December 8, 2020 from https://aperiodical.com/2020/06/conways-circle-theorem-a-proof-this-time-with-words

[3] Capitán, J. G. F. (2013). A generalization of the Conway circle, Forum Geom. 13: 191–195.

[4] Josefsson, M. (2011). More characterizations of tangential quadrilaterals, *Forum Geom*. 11: 65–82.

[5] Leung, A. and Lopez-Real, F. (2003). Properties of tangential and cyclic polygons: An application of circulant matrices, *International Journal of Mathematical Education* 34 (6): 859-870.

[6] Minculete, N. (2009). Characterizations of a tangential quadrilateral, *Forum Geom.* 9: 113-118.

[7] Mulcahy, C. Conway's Circle Theorem, retrieved December 8, 2020 from http://www.cardcolm.org/JHC.html#Circle

[8] Ryba, A. Theorem (Conway's Circle), retrieved December 8, 2020 from http://www.cardcolm.org/JHC.html